\documentclass[pdflatex,sn-mathphys-ay]{sn-jnl}

\usepackage{xurl}
\hypersetup{
  colorlinks=true,
  linkcolor=blue,
  citecolor=blue,
  urlcolor=blue
}
\usepackage{amsmath,amssymb}
\usepackage{graphicx}
\usepackage{xcolor}
\usepackage[normalem]{ulem}
\usepackage{dblfloatfix}
\usepackage[section]{placeins}

\title{Topology Optimization of Cooling Channels Using Dual-Type Moving Morphable Components}

\author*[1]{\fnm{Shunsuke} \sur{Hirotani}}\email{shunsuke.hirotani@nature-architects.com}
\author[1]{\fnm{Kunitaka} \sur{Shintani}}\email{kunitaka@nature-architects.com}
\author[1]{\fnm{Yoshikatsu} \sur{Furusawa}}\email{yoshikatsu.furusawa@nature-architects.com}
\author[2]{\fnm{Kentaro} \sur{Yaji}}\email{yaji@mech.eng.osaka-u.ac.jp}

\affil[1]{\orgname{Nature Architects Inc.}, \orgaddress{\street{3-8, Nihonbashi-Ningyocho, Chuo}, \city{Tokyo}, \postcode{103-0013}, \country{Japan}}}
\affil[2]{\orgname{The University of Osaka}, \orgaddress{\street{2-1, Yamadaoka}, \city{Suita}, \state{Osaka}, \postcode{565-0871}, \country{Japan}}}

\abstract{Efficient thermal management in high-power electronic devices requires cooling channel designs that provide high heat removal while satisfying strict spatial and manufacturing constraints. This study presents a two-stage hierarchical topology optimization framework for cooling channels based on the Moving Morphable Components (MMC) method. The optimization is performed sequentially: in the first stage, only wall components are optimized to establish the global flow network and insignificant components are removed; in the second stage, the global structure is fixed and fin components are optimized to improve local thermal performance. The method is coupled with a two-layer thermofluid model using the Brinkman approximation and solved with the adjoint sensitivity approach. Across multiple inlet pressure conditions, the proposed framework consistently generates designs with clear functional separation. The results demonstrate that exploring such clearly separated structures through a two-stage optimization strategy leads to a further reduction in the objective function. Compared with simultaneous MMC optimization and conventional density-based topology optimization, the proposed method produces geometries that are more interpretable, controllable, and suitable for manufacturing.}

\keywords{Topology optimization; Moving morphable components; Thermal-fluid problem; Cooling channels}

\begin{document}

\maketitle

\section{Introduction}
In recent years, the importance of thermal management has increased significantly due to continuous improvements in energy efficiency and the advancing performance of electronic devices \citep{Dhumal2023}. High-power batteries for electric vehicles (EVs), power electronic components, and highly integrated semiconductor devices are particularly prone to rapid temperature rises and hot-spot formation caused by locally intensified heat generation. Such temperature excursions can degrade performance, reduce reliability, and shorten operational lifetimes, potentially leading to catastrophic failure. Consequently, there is a strong demand for technologies that enable the design of efficient and reliable cooling structures under limited space and strict mechanical constraints. Among various cooling strategies, flow-channel-based liquid cooling systems are widely employed, and their geometric design has been recognized as a critical issue in next-generation electronic thermal management \citep{Garimella2008}. These studies highlight the importance of optimizing internal flow paths to enhance heat transfer while maintaining acceptable pressure losses. However, systematic design methodologies that can efficiently explore the trade-off between global flow control and local heat transfer enhancement remain an open challenge.

Over the past several decades, significant progress has been made in the size and shape optimization of fluid systems, exemplified by the seminal works \citep{Mohammadi2009}. However, these approaches are generally limited to optimizing pre-existing boundaries and do not allow for topological changes during the optimization process. To overcome this limitation, topology optimization has emerged as a powerful and flexible methodology \citep{Bendsoe1988,Bendsoe1989,Sethian2000}. By treating the material distribution within a prescribed domain as the design variable, topology optimization enables the automatic generation of structures that optimize objectives such as pressure loss, thermal performance, or structural compliance. Various formulations have been proposed, including homogenization-based, density-based, and level-set methods, and have been successfully applied to structural mechanics, heat conduction, fluid flow, and coupled thermofluid problems. \citet{Bendsoe2003} systematically established the theoretical and numerical foundations of topology optimization, particularly for density-based methods, thereby providing a basis for subsequent multiphysics extensions.

Topology optimization for fluid problems was initiated by \citet{Borrvall2003}, who addressed Stokes flow and introduced a framework for minimizing power dissipation under a fluid volume constraint. This framework was later extended to incompressible Navier--Stokes flows by \citet{Olesen2006}. Subsequent studies incorporated convective heat transfer and conjugate thermofluid effects, including flows at higher Reynolds numbers. For example, early efforts in topology optimization for heat transfer-related design include the work by \citet{Yoon2010}, who applied topology optimization to heat-dissipating structures with forced convective transport, demonstrating the importance of balancing conduction and convection. \citet{Alexandersen2015} subsequently investigated topology optimization for natural convection and showed the automatic generation of optimal heat sink geometries. \citet{Dilgen2018} extended topology optimization to practical high-Reynolds-number applications by combining turbulence modeling with adjoint sensitivity analysis. \citet{Yan2019} optimized solid--fluid bilayer and multilayer channel models while accounting for both thermal performance and pressure loss, demonstrating the emergence of complex and non-intuitive channel geometries. In addition, \citet{AlexandersenAndreasen2020} provided a comprehensive review of topology optimization methods for fluid flow and conjugate heat transfer problems, summarizing the state of the art and highlighting future research directions. These developments highlight the versatility and design freedom offered by topology optimization in thermal management.

Despite these advances, topology-optimized designs often insufficiently account for manufacturing constraints and practical industrial requirements. In density-based approaches, geometric smoothness and minimum feature sizes can be influenced using filtering and projection schemes; however, the explicit enforcement of constraints such as minimum channel widths, clearances, and curvature radii remains challenging. In fluid-related topology optimization, additional geometric requirements arise, including the need to ensure well-defined channel widths, suppress gray transition regions between solid and fluid, and maintain sufficient hydraulic diameters to avoid excessive pressure losses. Recent studies have also introduced geometry-control strategies tailored specifically to flow problems, such as density filters designed for pipe-like structures to regulate uniform wall thickness and cross-sectional consistency within fluid domains \citep{Choi2024}. To address geometric length-scale issues more broadly, \citet{Guest2004} proposed techniques for length-scale control, while \citet{Sigmund2007} introduced morphological filtering to obtain crisp and manufacturable designs. Robust formulations that account for manufacturing variability were developed by \citet{Wang2011} and \citet{Schevenels2011}. With the rapid advancement of additive manufacturing, orientation-dependent constraints such as overhang limitations have also been incorporated into topology optimization. \citet{Langelaar2017} proposed a filtering strategy tailored to additive manufacturing, and \citet{Allaire2017} established a theoretical framework for optimization under overhang constraints. Nevertheless, the complex geometries generated by topology optimization often require extensive post-processing, which can degrade performance and obscure the original design intent. This situation underscores the need for optimization frameworks that inherently restrict the design space to manufacturable and geometrically well-defined solutions.

The Moving Morphable Components (MMC) approach, proposed by \citet{Guo2014}, offers a promising alternative for integrating design intent with manufacturing constraints. Unlike conventional methods that implicitly represent material distribution, MMC employs explicitly defined geometric components whose positions, dimensions, and orientations are treated as design variables. As a result, constraints such as minimum channel widths or manufacturable orientations can be directly imposed through bounds on component parameters, providing clear and intuitive geometric controllability.

Several component-based topology optimization methods, including MMC, rely on feature-mapping techniques to project geometric information onto a fixed analysis mesh. These approaches can be broadly classified into combine-then-map and map-then-combine strategies. MMC belongs to the former category, in which individual topology description functions (TDFs) are first aggregated and subsequently mapped to material properties. In contrast, the geometry projection method proposed by \citet{Norato2015} maps each component individually prior to combination. A comprehensive comparison of these paradigms is provided by \citet{Wein2020}.

The MMC framework has been extended in several directions. \citet{Zhang2016} introduced tapered components with linearly varying thickness, while \citet{Guo2016} parameterized component centerlines using trigonometric functions to generate periodic wavy structures. \citet{Deng2016} incorporated joint-like constraints to explicitly model multi-link mechanisms. More recently, \citet{Hirotani2025} proposed a formulation based on cubic polynomial interpolation of neutral axes, in which the position vectors and cross-sectional orientation angles at both ends of each component are treated as design variables. This formulation ensures $C^1$ continuity between components and enables the representation of complex curved geometries. These developments demonstrate the high geometric expressiveness and flexibility of the MMC approach.

MMC has also been applied to thermofluid design. \citet{Yu2019} introduced an MMC-based framework in which fluid domains are explicitly represented and optimized, with channel centerlines described by polynomial functions. Their results showed that highly sinuous and non-intuitive cooling channels can be generated while maintaining explicit geometric control, thereby effectively bridging the gap between MMC and density-based approaches.

However, existing MMC-based thermofluid models exhibit notable limitations. In particular, polynomial-based channel representations are not well suited for modeling isolated solid features such as pin fins or baffles embedded within the flow. Because these formulations inherently assume continuous flow paths, they cannot easily represent closed solid contours. Moreover, when a single component type is used to simultaneously define global flow routing and local heat-transfer-promoting features, it becomes difficult to distinguish the roles of individual components in system-level flow control and local heat-transfer enhancement. This coupling complicates design interpretation, reduces robustness, and makes post-optimization refinement difficult. In addition, fine-scale features that also determine global flow paths can render system-level performance highly sensitive to local degradation or manufacturing tolerances.

To address these challenges, prior studies in topology optimization have explored hierarchical and multiscale design strategies, including staged optimization \citep{Rodrigues2002,Xia2014}, multiresolution methods \citep{Nguyen2010}, and two-scale optimization frameworks that explicitly account for manufacturability \citep{Schury2012}. Although these approaches are not specific to MMC formulations, they highlight the importance of separating global and local design roles.

Motivated by these considerations, the present study proposes a two-stage hierarchical MMC-based optimization framework that explicitly separates global flow-channel design from local heat-transfer enhancement. Unlike conventional staged or multiscale approaches, the proposed framework introduces two distinct component types within the MMC formulation: wall components, which define the global flow-channel skeleton, and fin components, which are embedded within the channels to locally enhance heat transfer. In the first stage, the global channel layout is optimized using wall components. In the second stage, with the global structure fixed, fin components are optimized to regulate local thermal performance. This clear separation of design scales enables the generation of robust and controllable cooling channel designs that inherently satisfy manufacturing constraints.

The remainder of this paper is organized as follows. Section~2 describes the thermofluid model and the proposed MMC-based optimization framework. Section~3 presents numerical examples under various design conditions and compares the results with existing methods. Finally, Section~4 concludes the paper and discusses future challenges and research directions.

\section{Optimization Framework}

\subsection{Proposed Optimization Framework}

\begin{figure*}[tbp]
  \centering
  \includegraphics[width=0.92\textwidth]{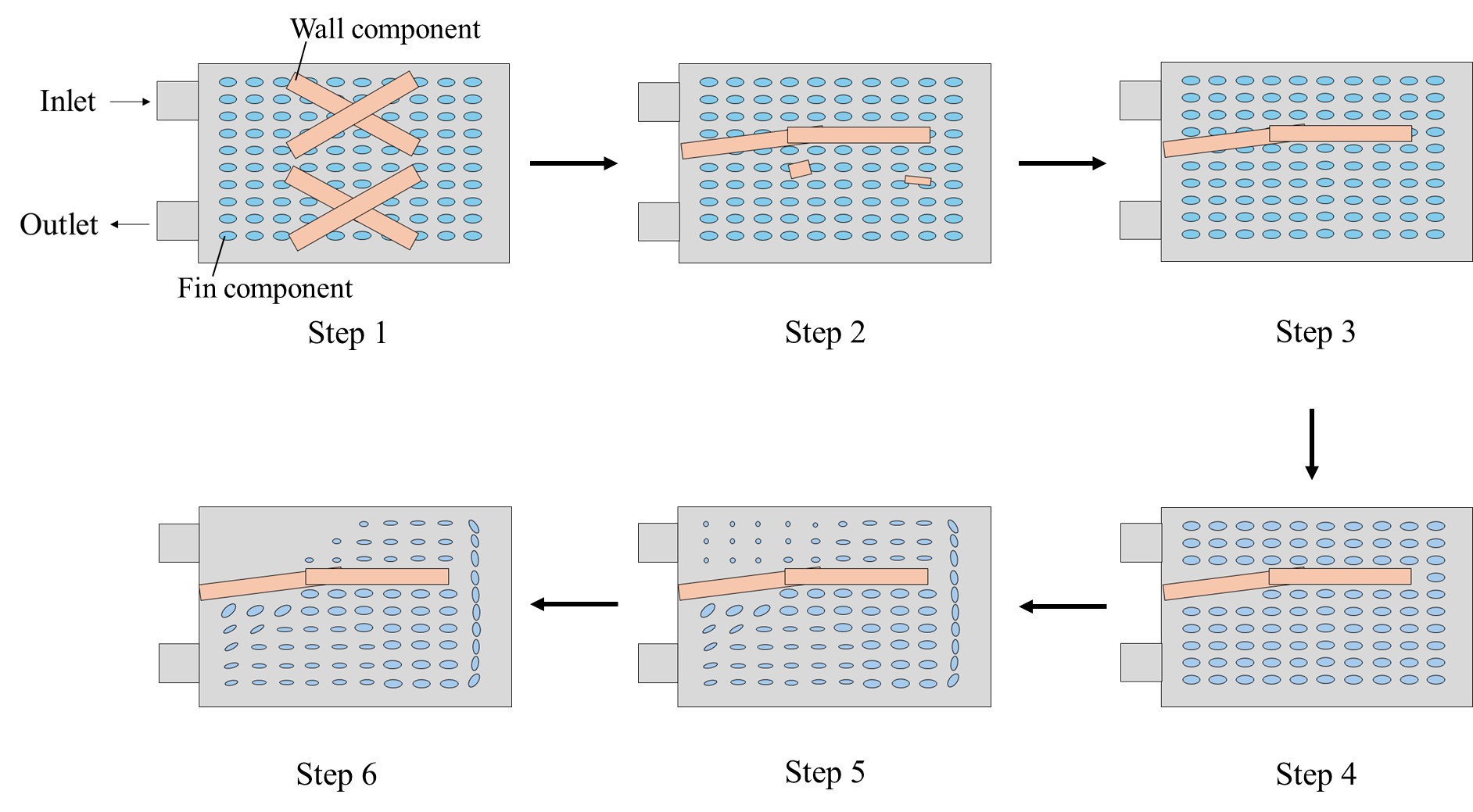}
  \caption{Conceptual diagram of the proposed optimization framework.}
  \label{fig:framework}
\end{figure*}

The proposed optimization framework is described. In this study, two types of MMC components with distinct functions are introduced: rectangular wall components that can freely move within the design domain, and elliptical fin components whose positions are fixed. Using these components, the optimization of the cooling channel is carried out through a two-stage hierarchical procedure. A conceptual diagram of the proposed optimization framework is shown in Fig.~\ref{fig:framework}.

The optimization is performed according to the following steps.

\noindent\textbf{Step 1.} Wall components and fin components are initially placed within the design domain.

\noindent\textbf{Step 2.} The global flow-channel network is optimized using only the wall components as design variables.

\noindent\textbf{Step 3.} Among the resulting wall components, those with sizes below a prescribed threshold that do not contribute to the formation of the global flow channels are removed.

\noindent\textbf{Step 4.} For the remaining wall components, fin components that are located closer than a predefined threshold distance are removed, leaving only candidate fin components for the second-stage optimization.

\noindent\textbf{Step 5.} The local heat transfer characteristics within the channels are optimized using only the fin components as design variables.

\noindent\textbf{Step 6.} Finally, fin components with sizes below a prescribed threshold are removed, and the resulting structure is adopted as the final design.

The primary advantages of separating the optimization into two stages can be summarized as follows:
\begin{enumerate}
\item[(i)] clear separation of functional roles across global and local spatial scales,
\item[(ii)] improved numerical stability and computational efficiency by progressively reducing the design space, and
\item[(iii)] generation of final designs that are physically interpretable and readily manufacturable.
\end{enumerate}

First, in the proposed framework, wall components are exclusively responsible for forming the global flow-channel network, whereas fin components are dedicated to controlling local heat transfer within the already established channels. By separating functions and spatial scales at the component level, the roles of design freedoms are clearly defined. As a result, if one wishes to modify global design features such as channel routing or branching locations, only the parameters of the wall components need to be adjusted. Conversely, if local heat transfer performance is to be tuned, only the dimensions or orientations of the fin components need to be modified. This separation enables localized and systematic post-optimization design modifications.

Second, the two-stage formulation improves numerical stability and efficiency by avoiding the simultaneous exploration of a strongly coupled, high-dimensional design space. In the first stage, the global channel network is formed using only wall components, and wall components below a size threshold that do not contribute to global channel formation are removed early in the process. Moreover, by removing fin candidates that are closer than a prescribed threshold to the finalized walls, the number of fin components considered in the second stage is further reduced. As a result, the second-stage optimization deals only with local structures under a fixed global layout, leading to a lower-dimensional and more weakly coupled optimization problem.

Third, the combination of two-stage optimization and threshold-based component removal improves both the physical interpretability and manufacturability of the final design. For both walls and fins, geometrically insignificant features below prescribed thresholds are systematically removed, eliminating nonessential microstructures. Consequently, the remaining wall components clearly define the channel skeleton, while the fin components are specialized for local heat transfer control within that skeleton.

\subsection{Thermofluid Model}

In this study, the two-layer thermofluid model is adopted, consisting of an upper thermofluid layer and a lower substrate layer \citep{Yan2019}. In the two-layer model, the upper layer treats both fluid and solid regions using a continuous density field $\gamma(x,y)\in[0,1]$ based on the Brinkman approximation, where the velocity in the solid region is suppressed through inverse permeability $\alpha(\gamma)$. The analysis is performed on a two-dimensional domain by solving the governing equations for fluid flow and heat transfer.

For the upper fluid layer, assuming incompressible flow, the continuity equation is given by
\begin{equation}
\nabla\cdot \mathbf{u} = 0,
\label{eq:continuity}
\end{equation}
where $\mathbf{u}$ denotes the velocity field. The momentum conservation equation based on the Brinkman approximation is expressed as
\begin{equation}
\frac{6}{7}\rho(\mathbf{u}\cdot\nabla)\mathbf{u}
= -\nabla p + \mu\nabla\cdot\left(\nabla\mathbf{u}+(\nabla\mathbf{u})^{T}\right) - \alpha(\gamma)\mathbf{u},
\label{eq:momentum}
\end{equation}
where $\rho$ is the fluid density, $p$ is the pressure and $\mu$ is the fluid dynamic viscosity.

The energy conservation equation for the upper layer is derived based on the two-layer model that assumes specific velocity and temperature profiles across the channel height, and is written as
\begin{equation}
\frac{2}{3}\rho C_{p}\mathbf{u}\cdot\nabla T_{t}
-\frac{49}{52}\nabla\cdot\bigl(k_{t}(\gamma)\nabla T_{t}\bigr)
-\frac{h(\gamma)}{2H_{t}(\gamma)}\left(T_{b}-T_{t}\right)=0,
\label{eq:energy_upper}
\end{equation}
where $C_{p}$ is the specific heat capacity of the fluid, $k_{t}(\gamma)$ is the effective thermal conductivity of the upper layer, $h(\gamma)$ is the interfacial heat transfer coefficient between the two layers, $H_t(\gamma)$ is the effective thickness of the upper layer, and $T_t$ and $T_b$ denote the representative temperatures of the upper and lower layers, respectively.

For the lower substrate layer, assuming a linear temperature distribution across the thickness direction, a two-dimensional heat conduction equation is formulated for the in-plane representative temperature $T_b$ as
\begin{equation}
\frac{k_b}{2}\nabla^2 T_b
+\frac{h(\gamma)}{2H_b}\left(T_b-T_t\right)
+\frac{q_0''}{2H_b}=0,
\label{eq:energy_lower}
\end{equation}
where $k_b$ is the thermal conductivity of the substrate material, $H_b$ is the half-thickness of the substrate, and $q_0''$ is the uniform heat flux applied from the bottom surface of the substrate. The second term represents the out-of-plane interfacial heat flux between the two layers, $q_{\mathrm{int}}=h(\gamma)\left(T_b-T_t\right)$. Through Eqs.~\eqref{eq:energy_upper} and \eqref{eq:energy_lower}, the convective--conductive heat transfer in the upper layer and the conductive heat transfer in the lower layer are thermally coupled, forming the two-layer thermofluid model. By solving the governing equations \eqref{eq:continuity}--\eqref{eq:energy_lower} together with the inlet and outlet pressure and temperature conditions and the velocity and temperature boundary conditions on the outer walls, the thermofluid fields of the cooling channel structures considered in this study are evaluated.

The material distribution $\gamma(x,y)$ is incorporated into these governing equations through effective material properties such as $\alpha(\gamma)$ and $k_t(\gamma)$. To smoothly and nonlinearly interpolate between the material properties of the fluid region ($\gamma=1$) and the solid region ($\gamma=0$), the Rational Approximation of Material Properties (RAMP) scheme is employed \citep{Alexandersen2023}. Specifically, the Brinkman resistance coefficient $\alpha(\gamma)$ is interpolated using a low resistance value $\alpha_f$ for the fluid phase and a high resistance value $\alpha_s$ for the solid phase as
\begin{equation}
\alpha(\gamma)=\alpha_f+(\alpha_s-\alpha_f)\frac{1-\gamma}{1+q_f\,\gamma},
\label{eq:ramp_alpha}
\end{equation}
where $q_f>0$ is a parameter controlling the sensitivity of the interpolation with respect to changes in the design variable. Similarly, the effective thermal conductivity of the upper layer $k_t(\gamma)$ is interpolated from the fluid thermal conductivity $k_f$ and the solid thermal conductivity $k_s$ as
\begin{equation}
k_t(\gamma)=k_f+(k_s-k_f)\frac{1-\gamma}{1+q_k\,\gamma},
\label{eq:ramp_k}
\end{equation}
where $q_k>0$ controls the slope of the interpolation.

\subsection{Formulation}

\subsubsection{Design Variables of Individual Components}

\begin{figure*}[tbp]
  \centering
  \includegraphics[width=0.92\textwidth]{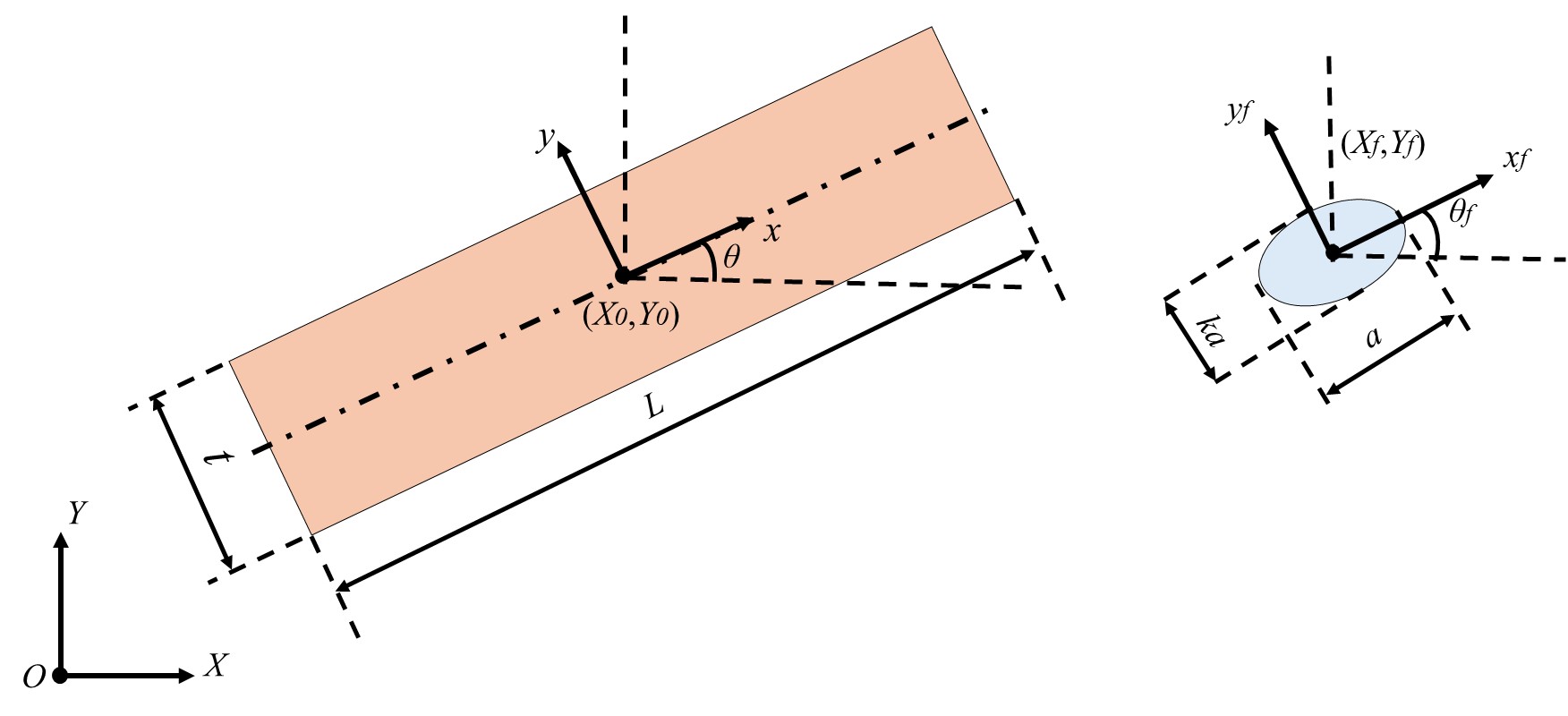}
  \caption{Wall components and fin components used in the optimization.}
  \label{fig:components}
\end{figure*}

In this study, two types of MMC components are introduced to represent the cooling channel geometry: rectangular wall components that can freely move within the design domain, and elliptical fin components whose positions are fixed, as illustrated in Fig.~\ref{fig:components}.

The wall components are introduced to form the global network and skeleton of the flow channels. Their design variables consist of the center coordinates $(X_{0}, Y_{0})$, the length in the longitudinal direction $L$, the thickness in the transverse direction $t$, and the rotation angle $\theta$.

In contrast, the fin components are introduced primarily to control the local heat transfer characteristics within the already established flow channels. In this study, the center coordinates $(X_{f}, Y_{f})$ of the fin components are fixed a priori. Consequently, the design variables of each fin component are defined as the semi-major axis length $a$, the ratio between the semi-minor and semi-major axes $0\le k\le 1$ (i.e., the semi-minor axis is given by $b = ka$), and the rotation angle $\psi$ of the ellipse.

\subsubsection{Topology Description Functions and Projection Functions}

The geometric shape of each MMC component is represented using a Topology Description Function (TDF), which varies smoothly with respect to the design variables \citep{Zhang2016}. To simplify the notation, let
\begin{align}
\Delta X_w &= X-X_0, & \Delta Y_w &= Y-Y_0, \label{eq:delta_w}\\
\Delta X_f &= X-X_f, & \Delta Y_f &= Y-Y_f. \label{eq:delta_f}
\end{align}
We further define the local coordinates of the wall component as
\begin{align}
s_w &= \Delta X_w \cos\theta + \Delta Y_w \sin\theta, \label{eq:sw}\\
n_w &= -\Delta X_w \sin\theta + \Delta Y_w \cos\theta, \label{eq:nw}
\end{align}
and those of the fin component as
\begin{align}
s_f &= \Delta X_f \cos\psi + \Delta Y_f \sin\psi, \label{eq:sf}\\
n_f &= -\Delta X_f \sin\psi + \Delta Y_f \cos\psi. \label{eq:nf}
\end{align}
Using these local coordinates, the TDF for a wall component is written as
\begin{equation}
\phi_{w}(X,Y)=1-\left(\frac{s_w}{L}\right)^{6}-\left(\frac{n_w}{t}\right)^{6},
\label{eq:tdf_wall}
\end{equation}
whereas the TDF for a fin component is written as
\begin{equation}
\phi_{f}(X,Y)=1-\left(\frac{s_f}{a}\right)^{2}-\left(\frac{n_f}{ka}\right)^{2}.
\label{eq:tdf_fin}
\end{equation}

In both cases, the TDF is constructed as a scalar field that takes positive values inside the component and negative values outside, with the zero level set representing the component boundary.

However, in the Brinkman-based unified field formulation adopted in the two-layer thermofluid model, the presence of fluid and solid phases is represented by a continuous density variable ranging between 0 and 1. In contrast, the TDF takes values in the range $(-\infty,+1)$ and therefore cannot be directly used as a material distribution variable. To address this issue, the TDF is projected onto an approximately binary function using a hyperbolic tangent mapping. Specifically, the projected TDFs are defined as
\begin{equation}
\tilde{\phi}_{w}(X,Y)=\frac{1}{2}\left(1-\tanh\left(\beta\phi_{w}(X,Y)\right)\right),
\label{eq:proj_wall}
\end{equation}
\begin{equation}
\tilde{\phi}_{f}(X,Y)=\frac{1}{2}\left(1-\tanh\left(\beta\phi_{f}(X,Y)\right)\right),
\label{eq:proj_fin}
\end{equation}
where $\beta$ is a parameter that controls the sharpness of the transition near the component boundaries.

Using these projected TDFs, the material distribution field $\gamma(X,Y)$ is defined as
\begin{equation}
\gamma(X,Y)=\tilde{\phi}_{w}(X,Y)\,\tilde{\phi}_{f}(X,Y).
\label{eq:gamma_def}
\end{equation}
The resulting material distribution $\gamma$ is then used in the interpolation of the Brinkman resistance term $\alpha(\gamma)$ and the effective thermo-physical properties in the governing equations.

\subsection{Optimization Formulation and Numerical Implementation}

\subsubsection{Optimization Problem}
In this study, the optimization objective is to minimize the normalized $p$-mean of the temperature in the solid layer (lower layer). Let $\Omega_b$ denote the lower-layer domain and $T_b(X,Y)$ the temperature field in the lower layer. The objective function $J$ is defined as
\begin{equation}
J=\left(
\frac{1}{|\Omega_b|}
\int_{\Omega_b}
T_b(X,Y)^{p}\,\mathrm{d}\Omega
\right)^{1/p}.
\label{eq:objective}
\end{equation}

Here, $p$ is an exponent used to approximate the maximum temperature, with $p\to\infty$ corresponding to $\max(T_b)$. In this study, $p=10$ is used \citep{Pejman2021}. As constraints, upper and lower bounds are imposed on the geometric parameters of the MMC components. The optimization is carried out according to the two-stage procedure described in Section~2.1.

\subsubsection{Numerical Implementation}
The following describes the numerical implementation details of the proposed optimization. The implementation is based on a finite element discretization of the two-layer thermofluid model together with adjoint sensitivities.

In the unified Brinkman-type formulation, the material distribution field $\gamma(\mathbf{x})$, defined from the projected topology description functions, is evaluated at the finite element level. The Brinkman resistance $\alpha(\gamma)$ and effective thermal conductivity $k_t(\gamma)$ are then computed via the RAMP interpolation. This element-wise representation allows the flow and heat-transfer equations to be assembled and solved on a fixed mesh without remeshing during the optimization.

The velocity, pressure, and temperature fields are discretized using first-order Lagrange finite elements. To ensure numerical stability, streamline-upwind/Petrov--Galerkin (SUPG) stabilization is applied to the governing equations in the thermofluid layer, including the Navier--Stokes and energy equations. The stationary thermofluid solution is computed at each design iteration, and adjoint problems associated with the objective functional are solved to evaluate gradients with respect to the component design variables.

Design updates are performed using the Method of Moving Asymptotes (MMA) \citep{Svanberg1987}. Each optimization iteration consists of a stationary state solve, adjoint sensitivity evaluation, and an MMA update of the design variables considered in each stage. The optimization process is terminated when the maximum normalized change in the design variables during one iteration becomes smaller than $10^{-3}$. In the present study, the workflow is implemented in COMSOL Multiphysics, using the stationary solver for the state analysis and the Optimization interface for adjoint-based sensitivities.

\section{Numerical Examples}

\subsection{Numerical Model}

\begin{figure}[tbp]
  \centering
  \includegraphics[width=\columnwidth]{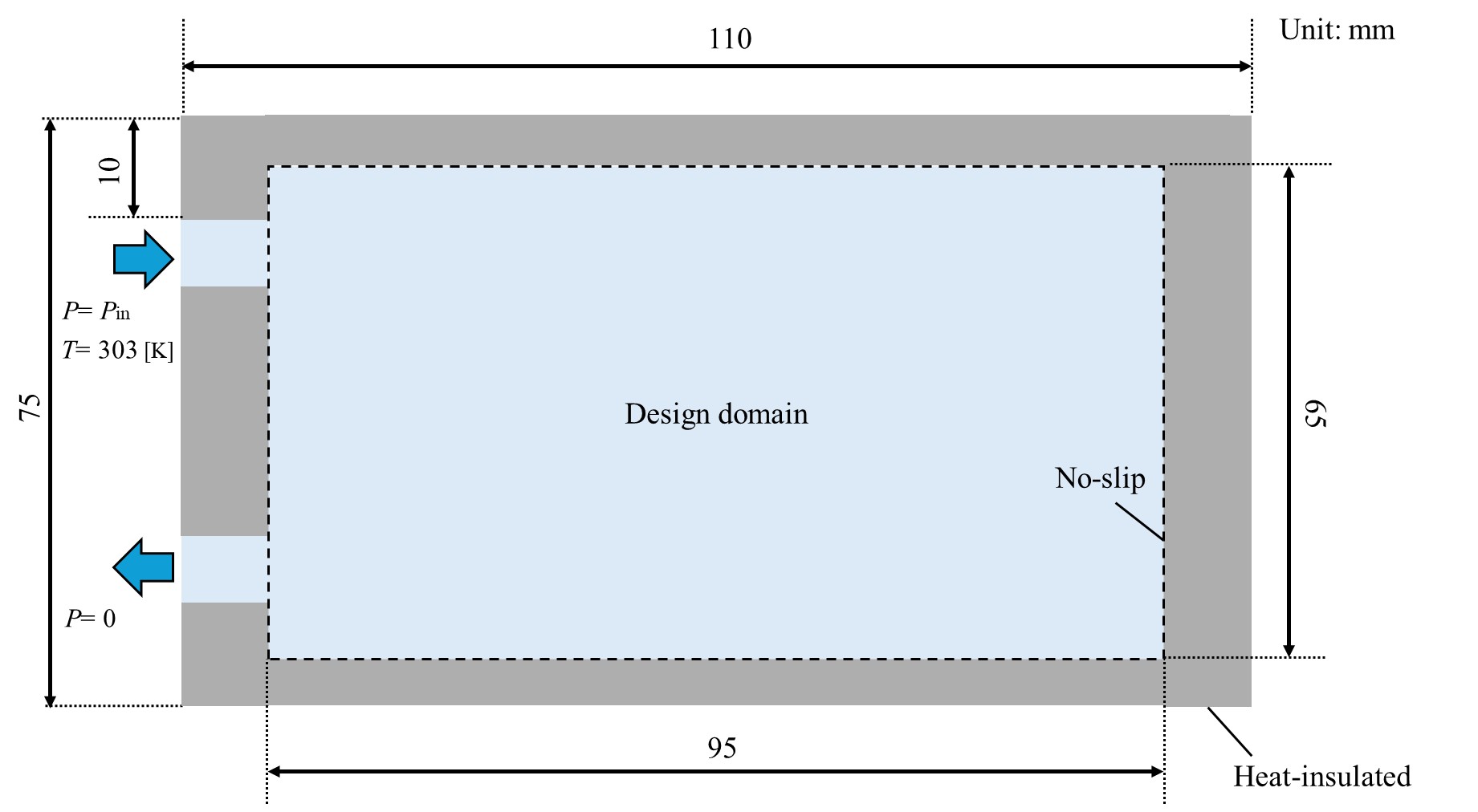}
  \caption{Geometry model and boundary conditions for the two-dimensional cooling device.}
  \label{fig:model_bc}
\end{figure}

In this section, numerical examples of cooling channel optimization using the proposed method are presented. The design domain is a rectangular region with an inlet located at the upper left and an outlet at the lower left, as shown in Fig.~\ref{fig:model_bc}. At the inlet, the pressure and temperature are prescribed as constant values ($p=p_{\mathrm{in}},\,T=303~\mathrm{K}$), while zero pressure is imposed at the outlet. No-slip conditions are applied on all wall boundaries, and adiabatic boundary conditions are imposed on the outer walls of the design domain. The model assumes a cooling plate attached to a heat-generating device, and a uniform heat flux of $q'' = 100~\mathrm{kW/m^2}$ is applied to the bottom surface of the substrate. Water and silicon are used as the fluid and solid materials, respectively. The physical properties employed in the simulations are summarized in Table~\ref{tab:material_properties}.

\begin{table}[tbp]
\centering
\caption{Material properties}
\begin{tabular}{lll}
\hline
Parameters & Values & Units \\
\hline
Water viscosity          & 0.001 & Pa $\cdot$ s \\
Water heat capacity      & 4180  & J/(kg $\cdot$ K) \\
Silicon heat capacity    & 942   & J/(kg $\cdot$ K) \\
Water conductivity       & 0.598 & W/(m $\cdot$ K) \\
Silicon conductivity     & 149   & W/(m $\cdot$ K) \\
\hline
\end{tabular}
\label{tab:material_properties}
\end{table}

\subsection{Results of the Two-Stage Optimization at Different Inlet Pressures}

\begin{figure}[tbp]
  \centering
  \includegraphics[width=0.6\textwidth]{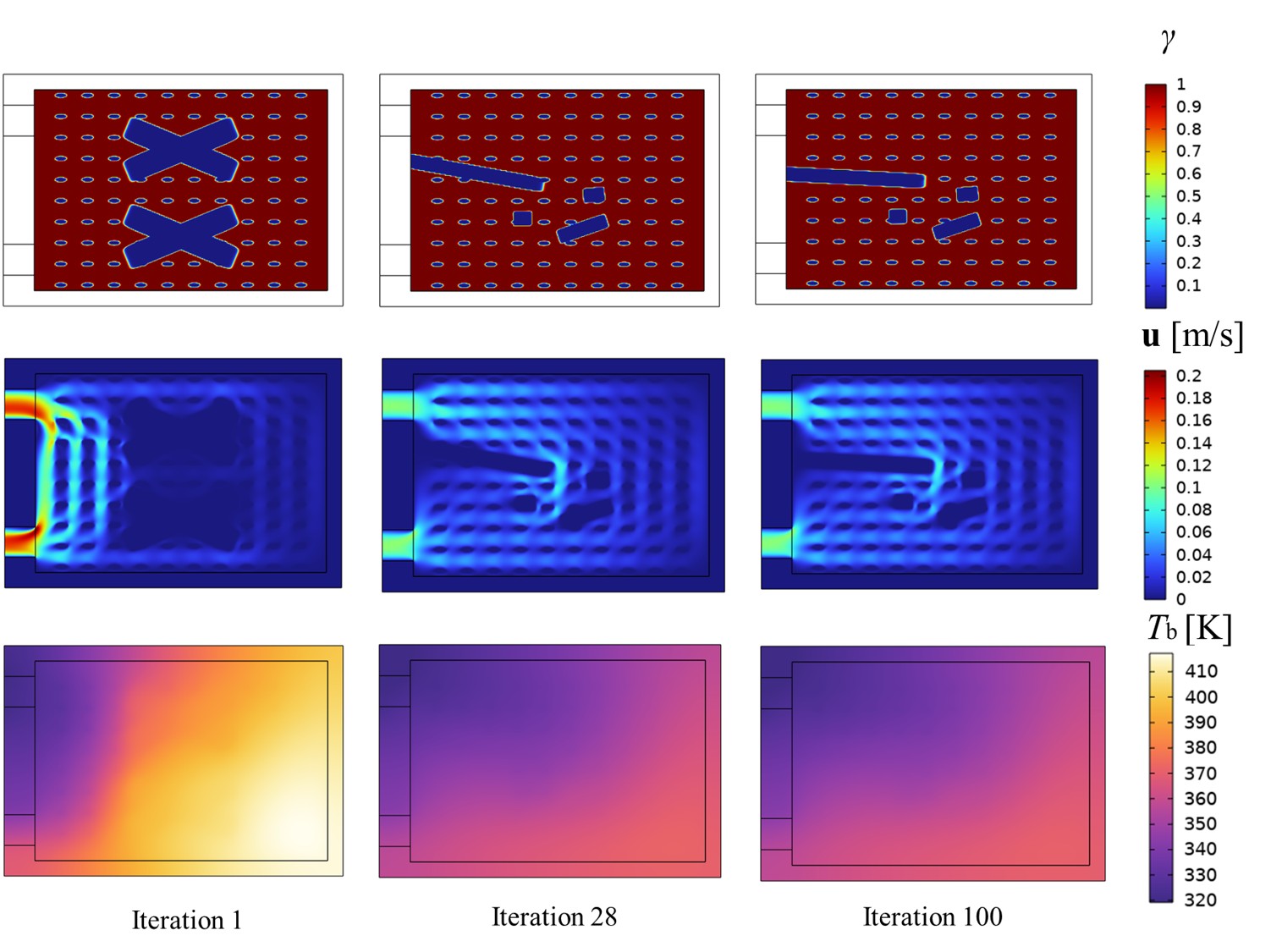}\\[-2mm]
  {\small (a)}\vspace{0.5mm}

  \includegraphics[width=0.6\textwidth]{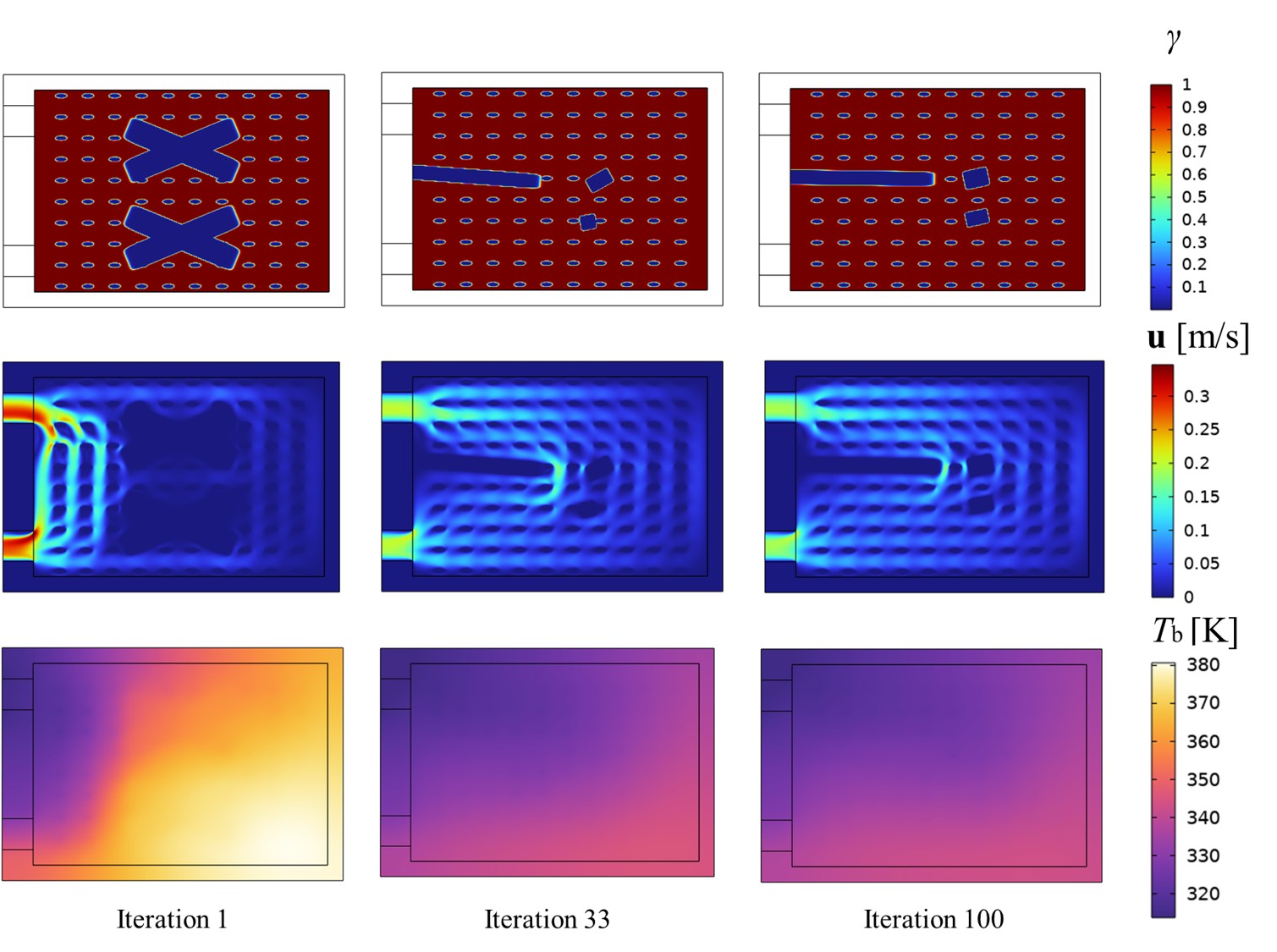}\\[-2mm]
  {\small (b)}\vspace{0.5mm}

  \includegraphics[width=0.6\textwidth]{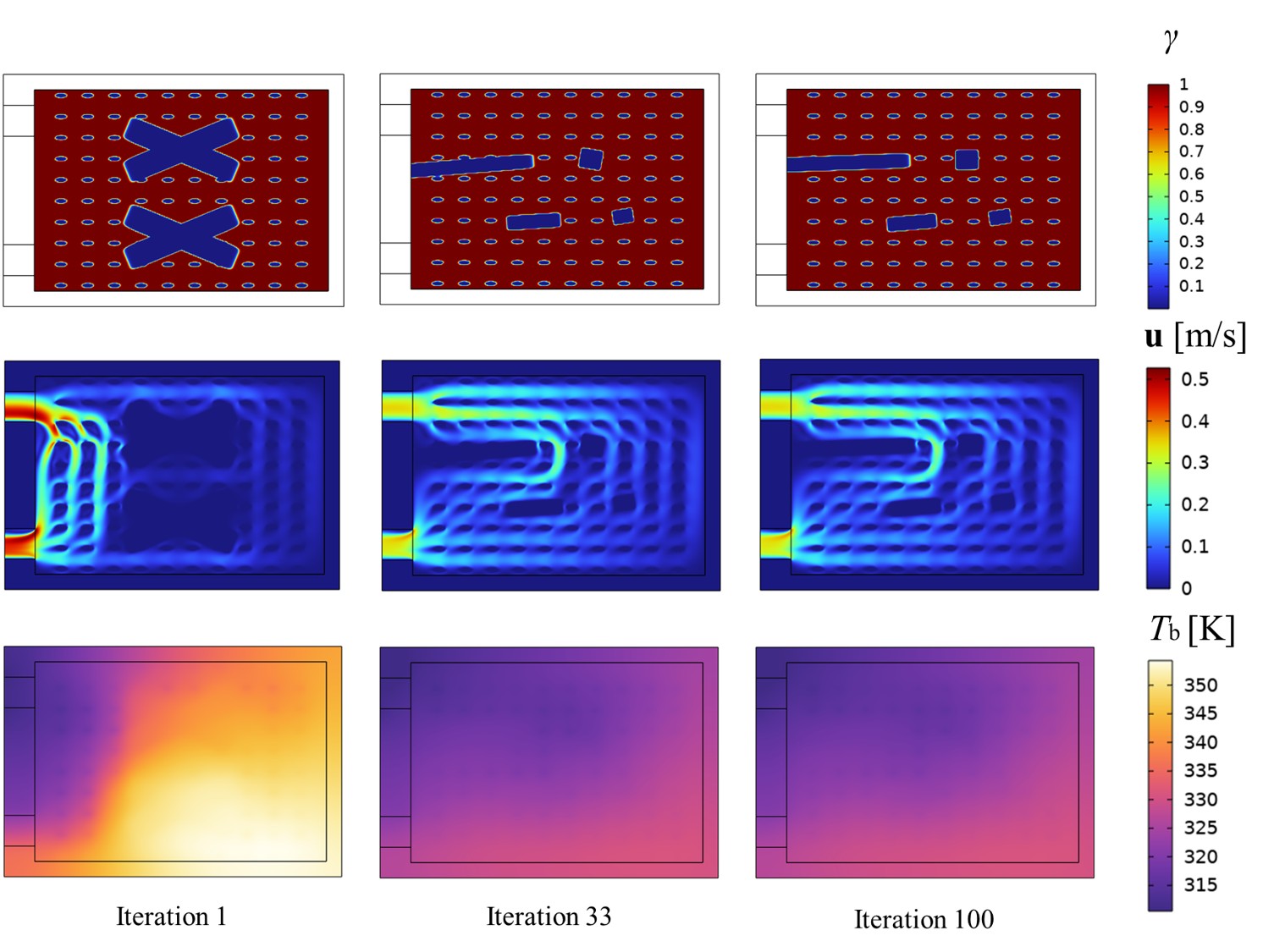}\\[-2mm]
  {\small (c)}
  \caption{Results at each inlet pressure during wall optimization: (a) $p_\text{in}=50$ Pa, (b) $p_\text{in}=100$ Pa, (c) $p_\text{in}=200$ Pa. The top row shows the material distribution, the middle row shows the flow velocity distribution, and the bottom row shows the temperature distribution of the solid layer.}
  \label{fig:wall_process}
\end{figure}

First, the optimization results for the wall components are presented for inlet pressures of $p_{\mathrm{in}}=50,\,100,$ and $200~\mathrm{Pa}$. The material distribution, velocity field, and solid-layer temperature distribution during the wall-optimization process are shown in Fig.~\ref{fig:wall_process} for (a) $p_{\mathrm{in}}=50~\mathrm{Pa}$, (b) $p_{\mathrm{in}}=100~\mathrm{Pa}$, and (c) $p_{\mathrm{in}}=200~\mathrm{Pa}$.

\begin{figure}[tbp]
  \includegraphics[width=\columnwidth]{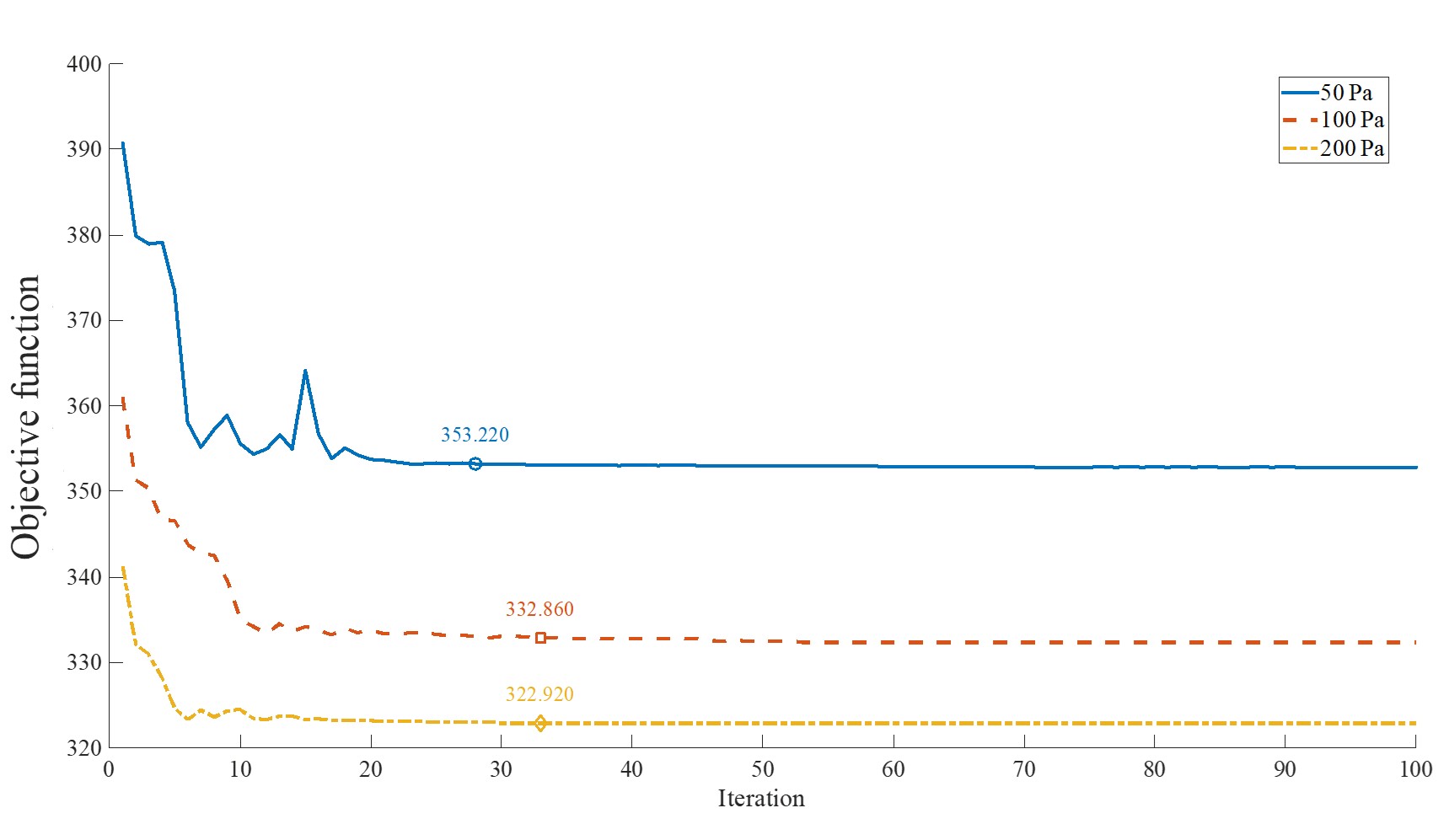}
  \caption{Objective function histories during wall-component optimization for multiple inlet pressure conditions.}
  \label{fig:wall_history}
\end{figure}

The corresponding histories of the objective function are shown in Fig.~\ref{fig:wall_history}. For all inlet pressures considered, the wall-component optimization suppresses the direct flow path from the inlet to the outlet and forms a channel that routes the flow around the right-hand side of the design domain, leading to a reduction in the objective function. Moreover, the solutions obtained after 100 iterations are almost identical to those obtained when the stopping criterion is satisfied, indicating that sufficient convergence is achieved before termination.

\begin{figure}[tbp]
  \centering
  \includegraphics[width=0.6\textwidth]{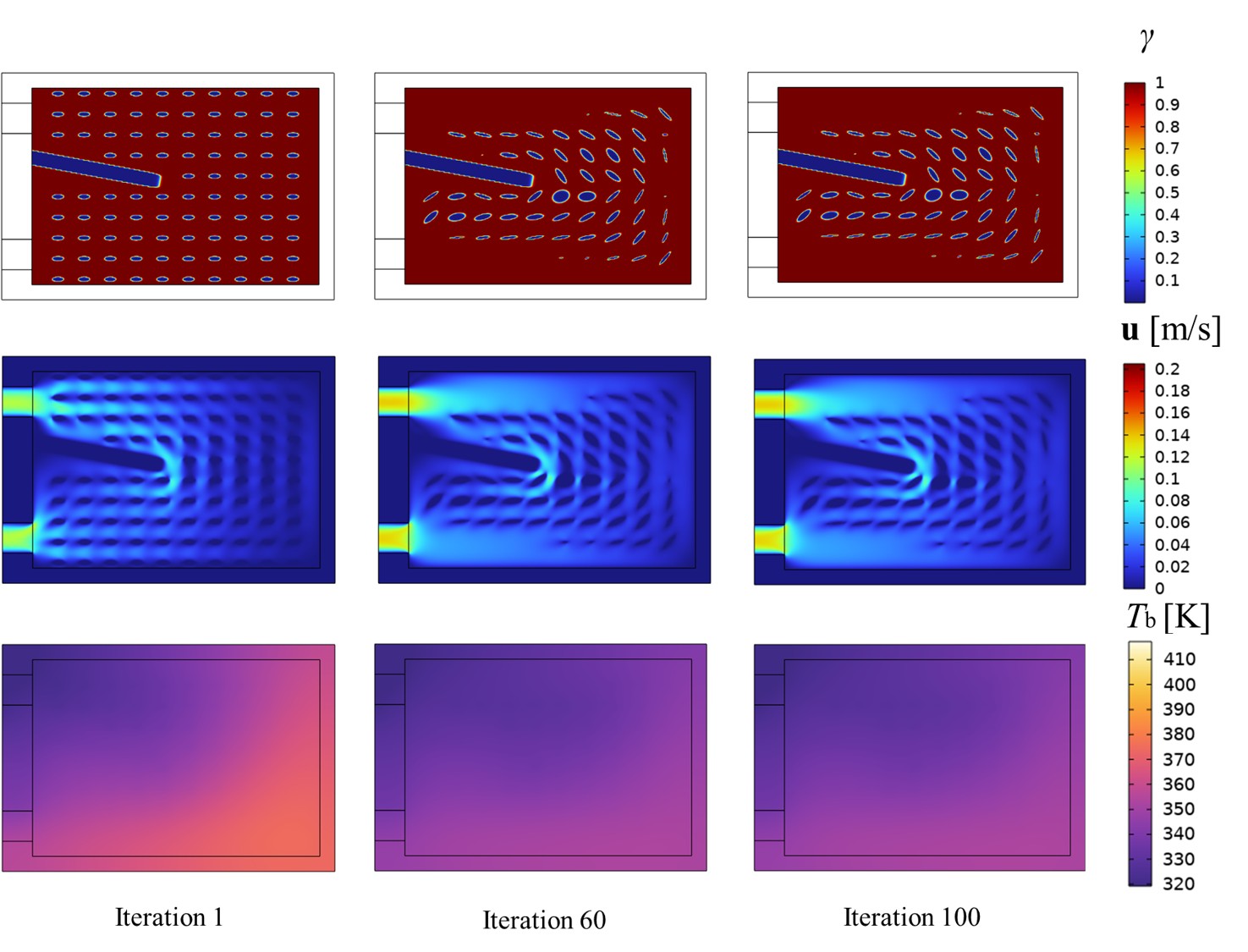}\\[-2mm]
  {\small (a)}\vspace{0.5mm}

  \includegraphics[width=0.6\textwidth]{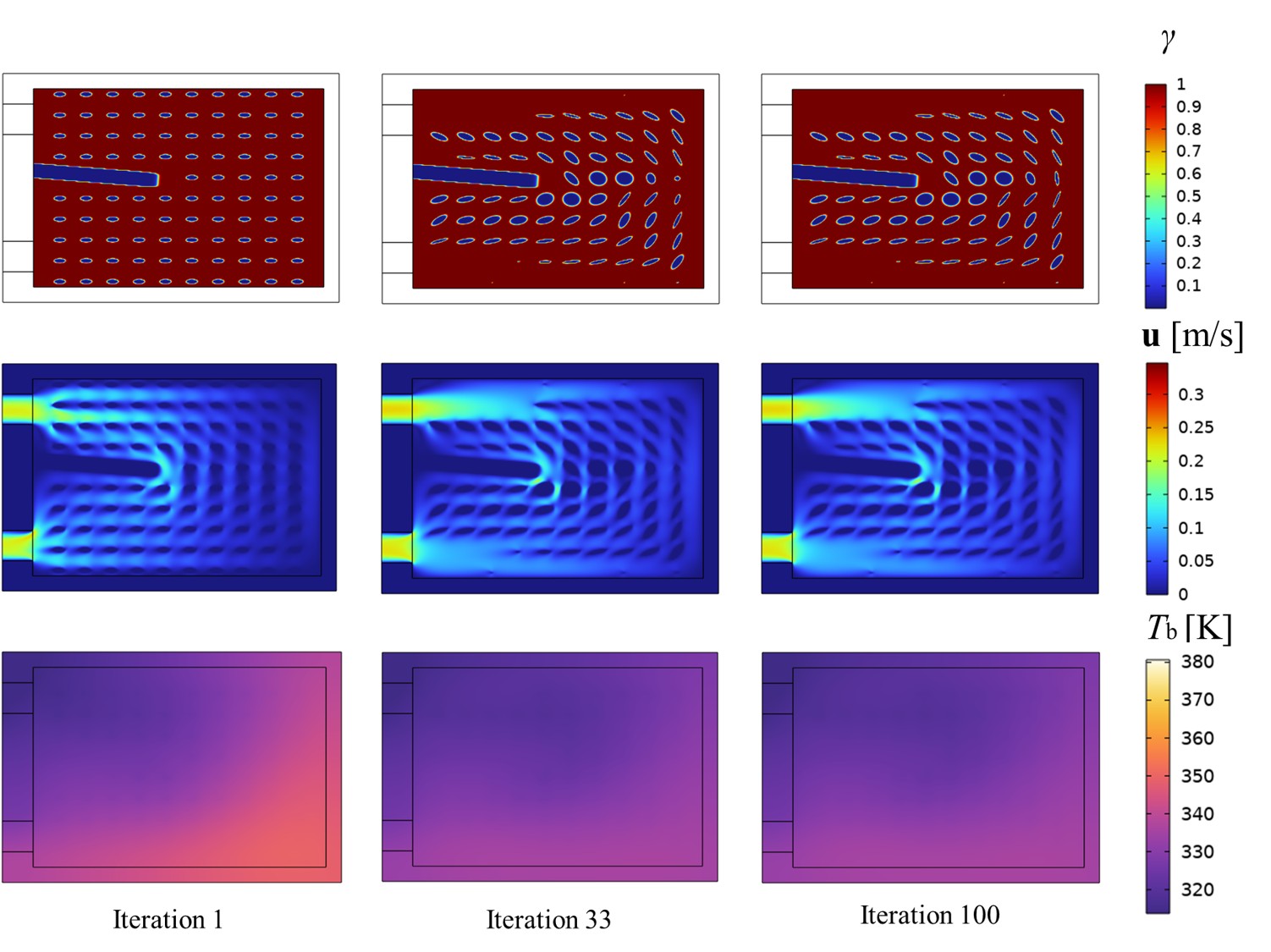}\\[-2mm]
  {\small (b)}\vspace{0.5mm}

  \includegraphics[width=0.6\textwidth]{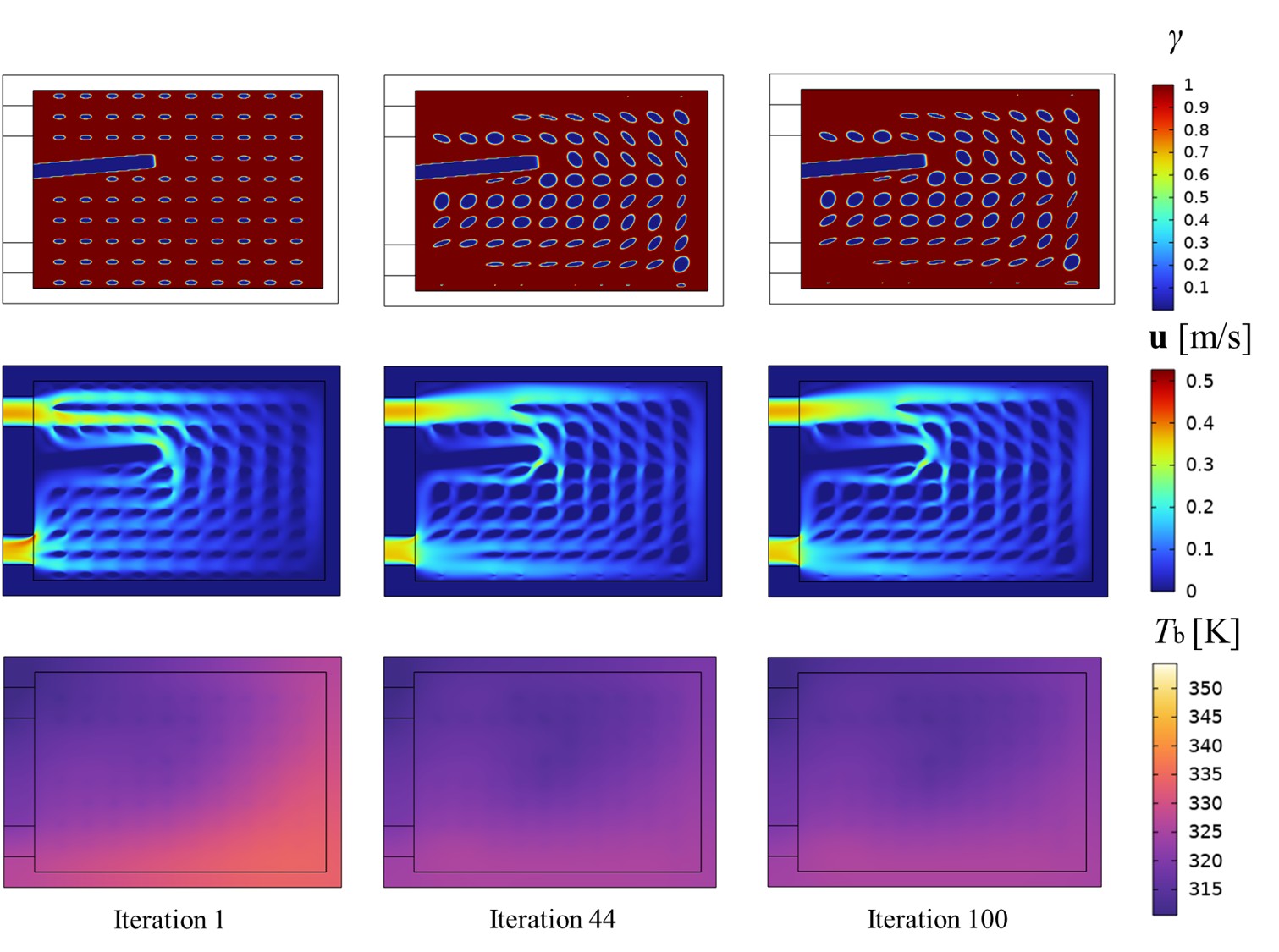}\\[-2mm]
  {\small (c)}
  \caption{Results at each inlet pressure during fin optimization: (a) $p_\text{in}=50$ Pa, (b) $p_\text{in}=100$ Pa, (c) $p_\text{in}=200$ Pa. The top row shows the material distribution, the middle row shows the flow velocity distribution, and the bottom row shows the temperature distribution of the solid layer.}
  \label{fig:fin_process}
\end{figure}

Next, the optimization results for the fin components are presented. Figure~\ref{fig:fin_process} shows the material distribution, velocity field, and solid-layer temperature distribution during fin optimization for (a) $p_{\mathrm{in}}=50~\mathrm{Pa}$, (b) $p_{\mathrm{in}}=100~\mathrm{Pa}$, and (c) $p_{\mathrm{in}}=200~\mathrm{Pa}$.

\begin{figure}[tbp]
  \centering
  \includegraphics[width=\columnwidth]{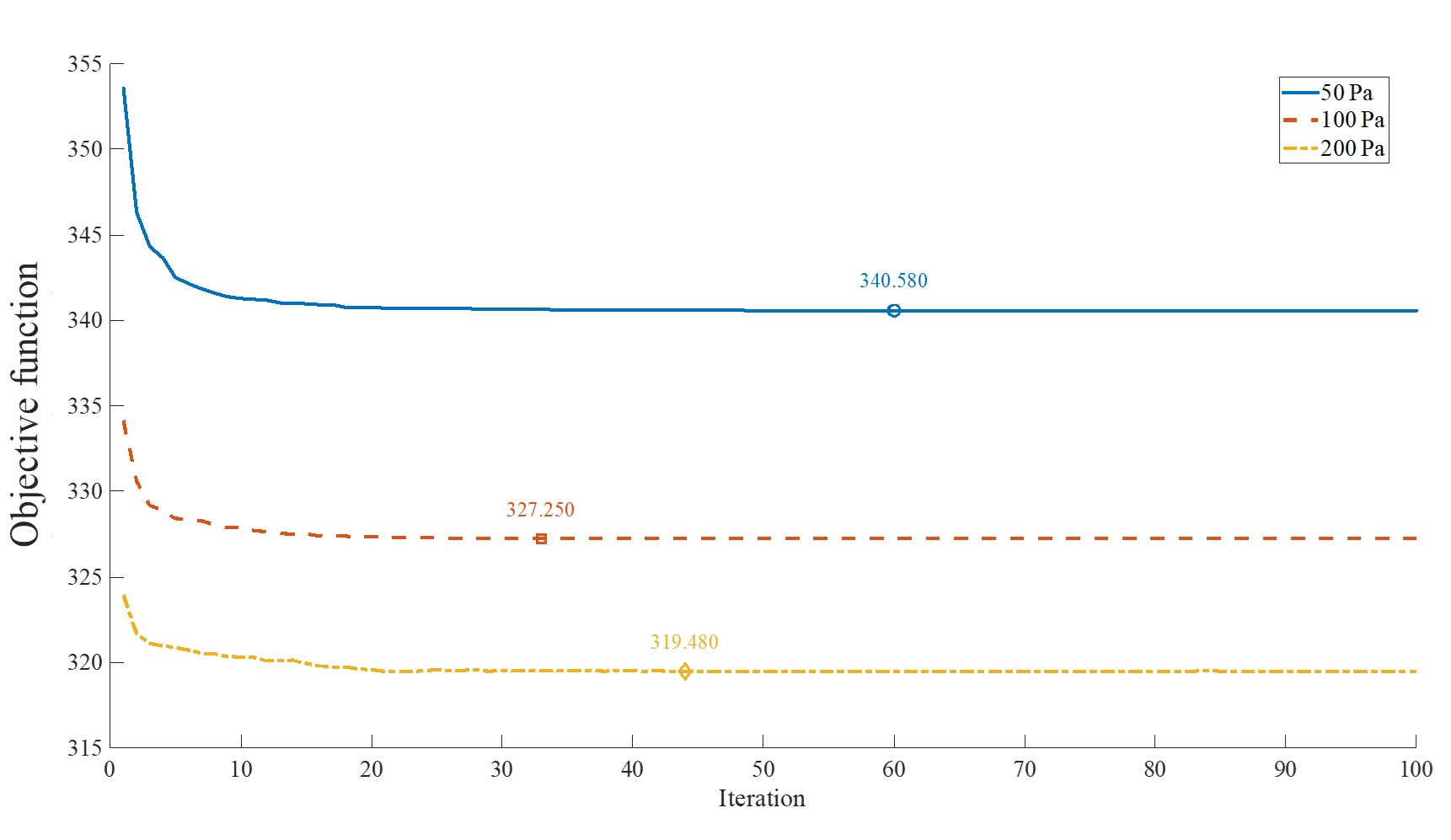}
  \caption{Objective function histories during fin-component optimization for multiple inlet pressure conditions.}
  \label{fig:fin_history}
\end{figure}

Figure~\ref{fig:fin_history} shows the corresponding histories of the objective function. The initial design for fin optimization is obtained by removing small wall components from the wall-optimization result. In this study, wall components with lengths smaller than 20\% of the longitudinal dimension of the design domain are removed. At all inlet pressures, the fin components effectively regulate local flow patterns, yielding further reductions in the objective function. It is also observed that fins tend to be more sparsely distributed as the inlet pressure decreases. This trend can be interpreted as follows. At higher inlet pressures, convection is stronger and the flow velocities are higher; therefore, flow curvature and recirculation induced by fins enhance mixing and promote thermal-boundary-layer thinning, which increases heat transfer. At lower inlet pressures, the flow becomes relatively weaker and more diffusion-dominated; in this regime, the convective contribution is limited while the additional flow resistance introduced by fins becomes more influential, making dense fin distributions less beneficial.

\begin{figure}[tbp]
  \centering
  \includegraphics[width=\columnwidth]{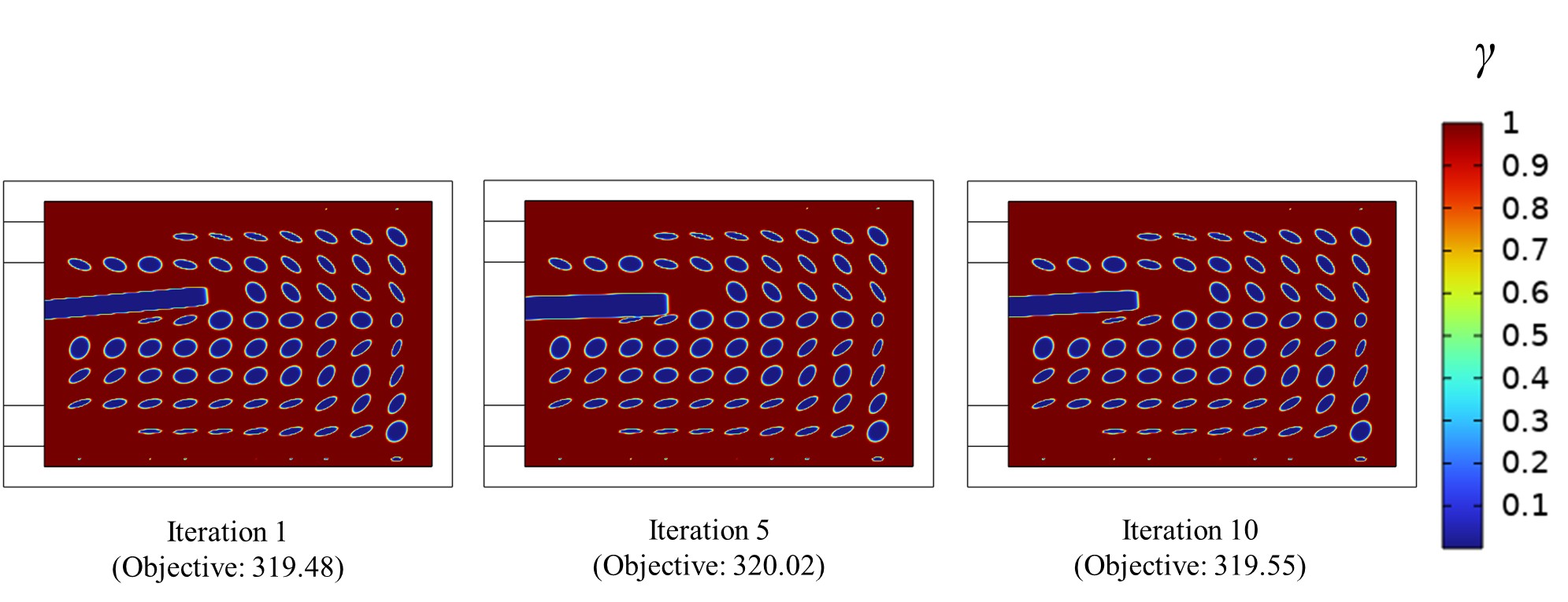}
  \caption{Optimization process when wall components are re-optimized.}
  \label{fig:reopt_wall}
\end{figure}

To further assess convergence and coupling between stages, the fin-optimization result is fed back into the wall optimization. Figure~\ref{fig:reopt_wall} shows the material distribution obtained when the wall components are optimized again at an inlet pressure of $200~\mathrm{Pa}$. The material distribution remains almost unchanged, indicating that the wall-component optimum is largely insensitive to the fin-optimization outcome and that the two-stage procedure is stable.

\subsection{Comparison Between Simultaneous Optimization and Two-Stage Optimization}

\begin{figure}[tbp]
  \centering
  \includegraphics[width=\columnwidth]{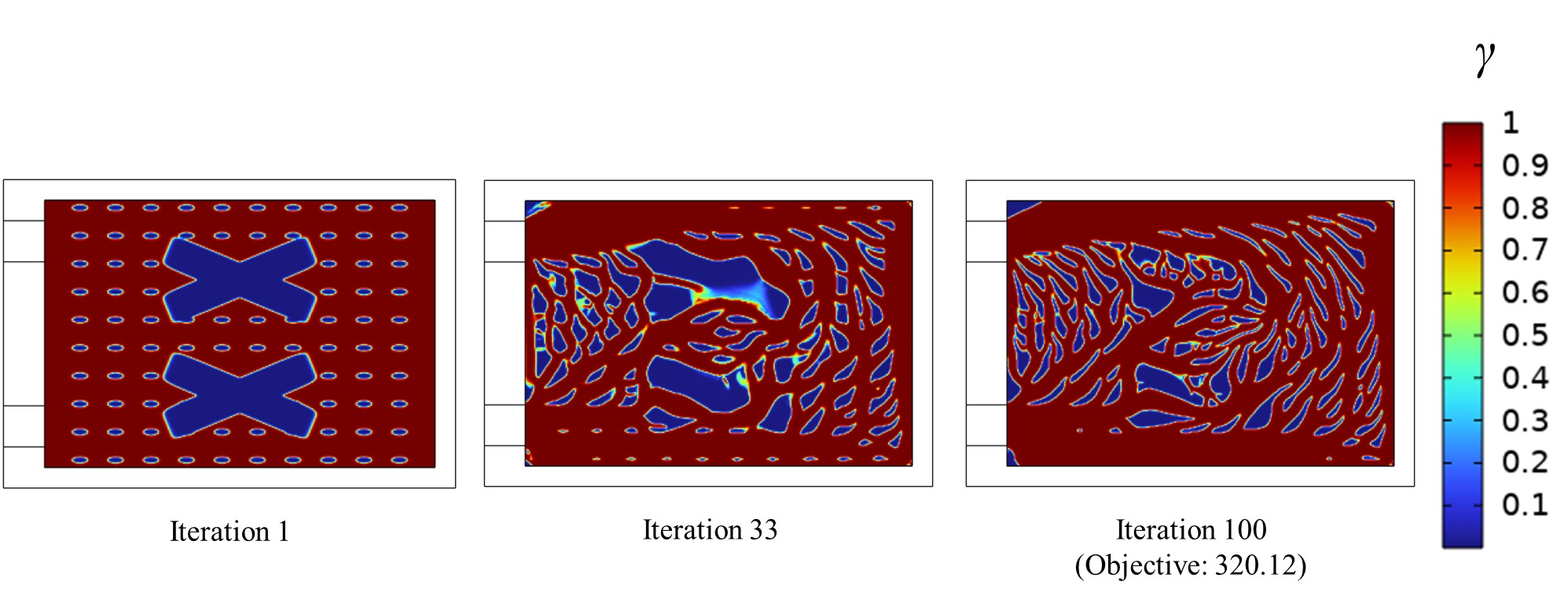}
  \caption{Simultaneous optimization process of wall and fin components.}
  \label{fig:simul_process}
\end{figure}

Next, the results of simultaneous optimization, in which wall and fin components are optimized concurrently, are compared with those obtained using the proposed two-stage optimization framework. The result of the simultaneous optimization at an inlet pressure of $200~\mathrm{Pa}$ is shown in Fig.~\ref{fig:simul_process}.

\begin{figure}[tbp]
  \centering
  \includegraphics[width=\columnwidth]{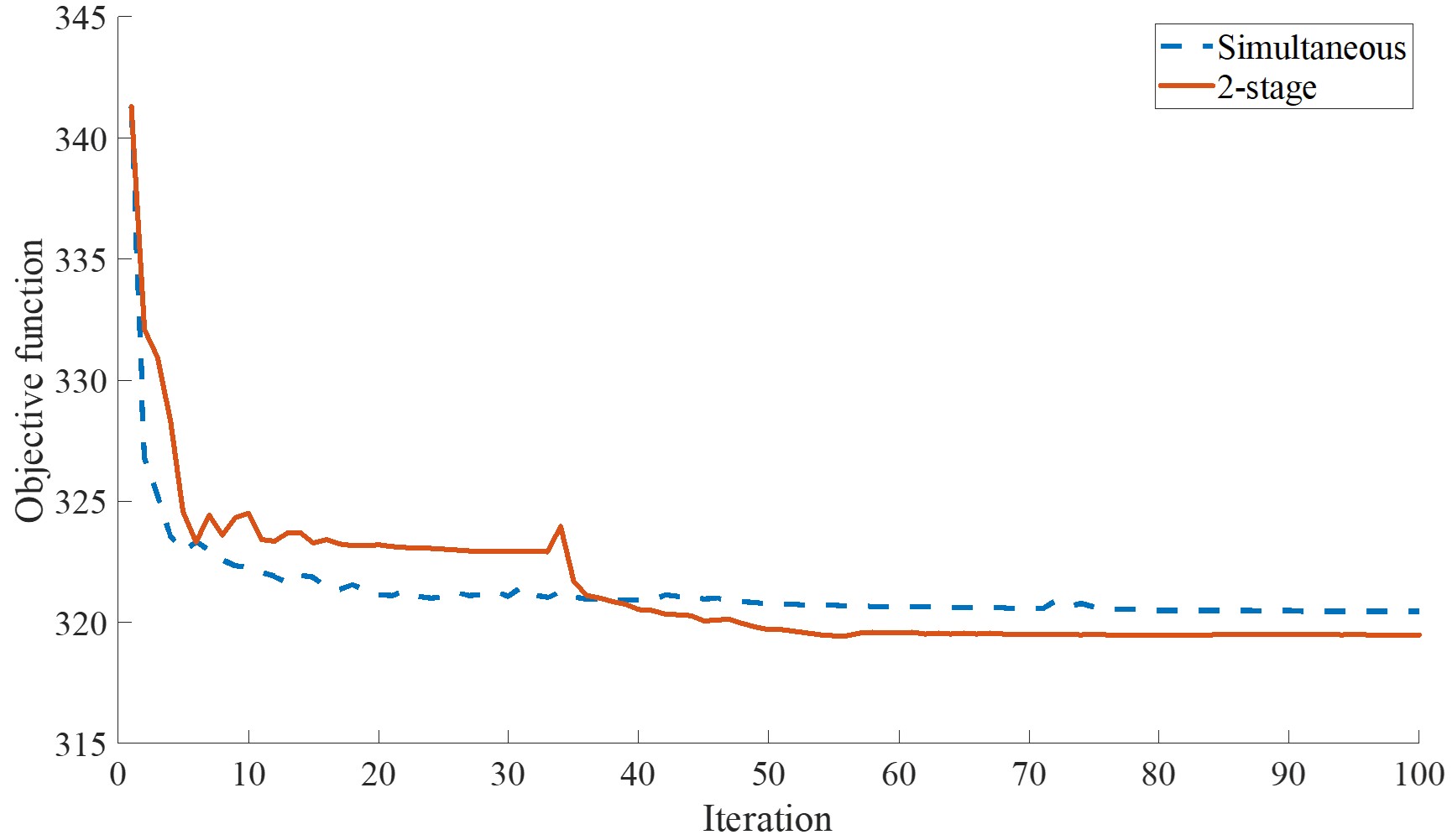}
  \caption{Comparison of objective function histories between the two-stage optimization and the simultaneous optimization.}
  \label{fig:simul_history}
\end{figure}

From the material distribution in Fig.~\ref{fig:simul_process}, it can be observed that the two component types no longer maintain a clear division of functional roles in the simultaneous optimization. Specifically, wall components are partly used for local heat-transfer enhancement (originally the role of fins), while fin components contribute to forming the global flow-channel network. Consequently, the separation between global channel formation and local heat-transfer enhancement becomes unclear. In contrast, the proposed method determines the global flow-channel skeleton using only wall components in the first stage, and then optimizes only fin components within the established channels in the second stage. This staged process preserves clear functional separation and progressively reduces the design space. A comparison of the objective-function histories is presented in Fig.~\ref{fig:simul_history}. From Fig.~\ref{fig:simul_history}, the proposed two-stage optimization achieves a lower final objective-function value. Although the simultaneous optimization reduces the objective function more rapidly in the early iterations because both wall and fin variables are updated simultaneously, it eventually converges to a higher local minimum than the two-stage approach.

\subsection{Comparison with Density-Based Topology Optimization}

\begin{figure}[tbp]
  \centering
  \includegraphics[width=\columnwidth]{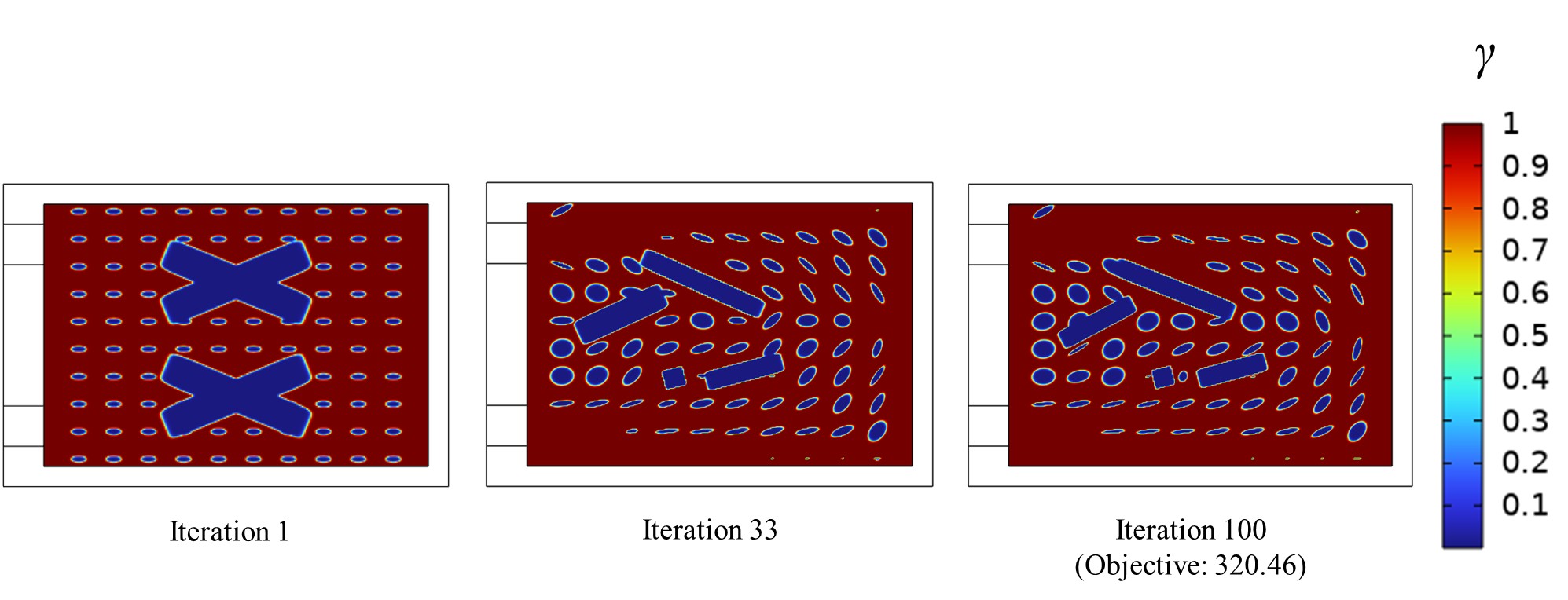}
  \caption{Result of density-based topology optimization.}
  \label{fig:density_result}
\end{figure}

\begin{figure}[tbp]
  \centering
  \includegraphics[width=\columnwidth]{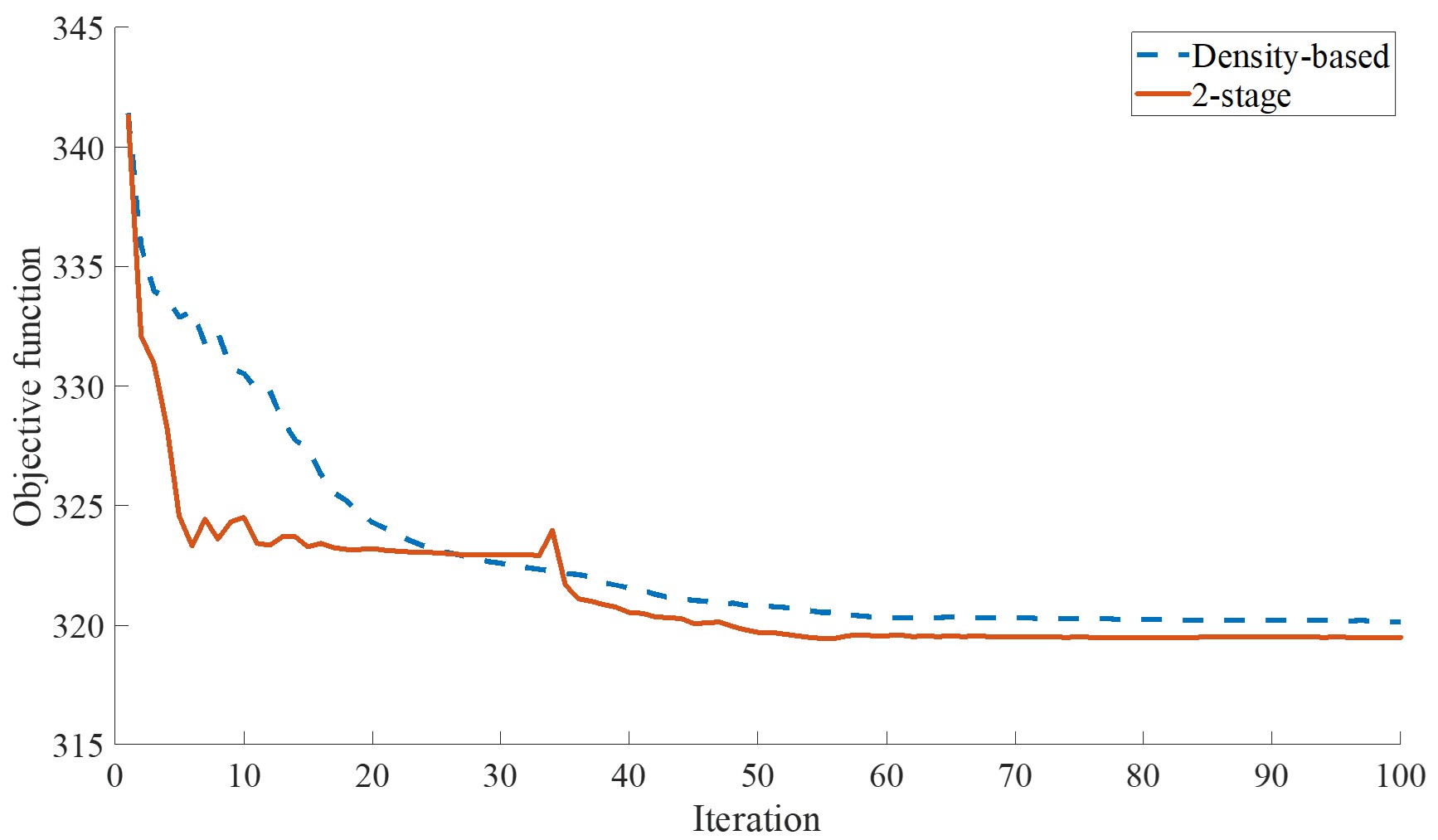}
  \caption{Comparison of objective function histories between the two-stage optimization and density-based topology optimization.}
  \label{fig:density_history}
\end{figure}

Finally, the proposed two-stage optimization framework is compared with conventional density-based topology optimization. The result of the density-based topology optimization at an inlet pressure of $200~\mathrm{Pa}$ is shown in Fig.~\ref{fig:density_result}, whereas Fig.~\ref{fig:density_history} presents the comparison of the corresponding objective-function histories.

From the material distribution in Fig.~\ref{fig:density_result}, it can be observed that density-based topology optimization treats the design as a pixel-wise continuous field. This provides extremely high geometric freedom and often yields complex and non-intuitive channel geometries. In contrast, the proposed MMC-based approach represents the geometry using explicit component parameters, resulting in designs with clearer physical interpretation and higher compatibility with CAD-based redesign and reuse.

From Fig.~\ref{fig:density_history}, the proposed two-stage optimization converges more rapidly and achieves a lower final objective-function value. These results demonstrate that the proposed method can attain performance comparable to or better than density-based topology optimization while maintaining superior geometric controllability and practical manufacturability.

\FloatBarrier

\section{Conclusion}
In this study, a hierarchical MMC-based topology optimization framework was proposed for the design of cooling channel systems, in which wall and fin components are optimized in two successive stages. Rectangular wall components are first introduced to establish the global skeleton of the flow-channel network, and elliptical fin components are subsequently optimized within this predefined skeleton to regulate local heat transfer. By explicitly separating global flow-channel formation and local heat transfer enhancement across different spatial scales, the proposed strategy achieves a structured decomposition of the design space and clarifies the physical roles of individual components.

Numerical investigations under inlet pressures of 50, 100, and 200 Pa demonstrate that the proposed approach consistently generates cooling channel configurations with a clear functional allocation: wall components predominantly determine the global flow topology, while fin components enhance local mixing and heat transfer. It is further observed that the optimal fin density decreases as the inlet pressure is reduced. This tendency reflects the shift in the trade-off between pressure loss and heat transfer enhancement: under low-flow conditions, the additional flow resistance induced by fins becomes comparatively more influential, leading the optimization process to suppress excessive fin insertion. These results confirm that the proposed method captures physically consistent design tendencies rather than merely producing numerical artifacts.

Compared with density-based topology optimization, the two-stage MMC framework involves a lower degree of geometric freedom due to its explicit component parameterization. Nevertheless, the resulting channel geometries exhibit simpler and more interpretable structural patterns that directly correspond to meaningful geometric parameters. From an optimization perspective, the hierarchical separation effectively regularizes the design space by reducing strong cross-couplings between global and local design variables. Despite the reduced geometric freedom, the objective function defined by the $p$-mean temperature of the solid layer decreases more rapidly during optimization and converges to lower final values than those obtained by density-based methods. This indicates that an appropriately structured representation can enhance both numerical stability and physical effectiveness without relying on excessive geometric complexity.

In contrast, when wall and fin components are optimized simultaneously, their structural roles tend to overlap and compete within a strongly coupled design space. Such coupling obscures the distinction between global flow formation and local heat transfer mechanisms, and leads to less effective objective reduction. The results therefore highlight the importance of hierarchically decomposing multi-scale design problems to improve sensitivity conditioning and to balance competing physical effects in thermofluid optimization.

A central advantage of the MMC framework lies in its flexibility to tailor component shape descriptions and design-variable relationships according to specific functional and manufacturing requirements. Although the present study focuses on wall and fin components, the methodology can be systematically extended to incorporate additional component types that enforce geometric constraints, manufacturability considerations, or prescribed connectivity conditions.

Future work will extend the proposed approach to higher Reynolds number regimes, transient thermal loading, and fully three-dimensional formulations. Furthermore, the direct embedding of realistic manufacturing constraints into MMC component definitions, as well as extensions to multi-objective and robust optimization settings, will further enhance its practical applicability. Through these developments, hierarchical component-based topology optimization has the potential to provide a physically interpretable and practically implementable framework for advanced cooling channel design.

\section*{Declarations}

\subsection*{Funding}
The authors received no external funding for this study.

\subsection*{Conflict of interest}
On behalf of all authors, the corresponding author states that there is no conflict of interest.

\subsection*{Author contributions}
Shunsuke Hirotani: Conceptualization, Methodology, Software, Formal analysis, Investigation, Visualization, Writing--original draft.
Kunitaka Shintani: Methodology, Supervision, Writing--review and editing.
Yoshikatsu Furusawa: Validation, Writing--review and editing.
Kentaro Yaji: Methodology, Supervision, Writing--review and editing.

\subsection*{Data availability}
The data generated and/or analyzed during the current study are available from the corresponding author on reasonable request.

\subsection*{Code availability}
The code and computational model used in the current study, including the COMSOL Multiphysics setup and optimization workflow, are available from the corresponding author on reasonable request.

\subsection*{Replication of results}
The results presented in this paper can be reproduced based on the numerical model, optimization formulation, and implementation details described in Sections 2 and 3. In particular, the governing equations, material interpolation scheme, MMC parameterization, optimization procedure, and boundary conditions are provided to support reproducibility. The data, code, and computational model used in this study are available from the corresponding author on reasonable request.


\begin{thebibliography}{99}

\bibitem[Alexandersen et~al.(2015)]{Alexandersen2015}
Alexandersen J, Aage N, Andreasen CS, Sigmund O (2015) Topology optimisation for natural convection problems. Int J Numer Methods Fluids 76(10):699--721. \url{https://doi.org/10.1002/fld.3954}

\bibitem[Alexandersen and Andreasen(2020)]{AlexandersenAndreasen2020}
Alexandersen J, Andreasen CS (2020) A review of topology optimisation for fluid-based problems. Fluids 5(1):29. \url{https://doi.org/10.3390/fluids5010029}

\bibitem[Alexandersen(2023)]{Alexandersen2023}
Alexandersen J (2023) A detailed introduction to density-based topology optimisation of fluid flow problems with implementation in MATLAB. Struct Multidisc Optim 66(1):12. \url{https://doi.org/10.1007/s00158-022-03420-9}

\bibitem[Allaire et~al.(2017)]{Allaire2017}
Allaire G, Dapogny C, Estevez R, Faure A, Michailidis G (2017) Structural optimization under overhang constraints imposed by additive manufacturing technologies. J Comput Phys 351:295--328. \url{https://doi.org/10.1016/j.jcp.2017.09.041}

\bibitem[Bends{\o}e(1989)]{Bendsoe1989}
Bends{\o}e MP (1989) Optimal shape design as a material distribution problem. Struct Optim 1(4):193--202. \url{https://doi.org/10.1007/BF01650949}

\bibitem[Bends{\o}e and Kikuchi(1988)]{Bendsoe1988}
Bends{\o}e MP, Kikuchi N (1988) Generating optimal topologies in structural design using a homogenization method. Comput Methods Appl Mech Eng 71(2):197--224. \url{https://doi.org/10.1016/0045-7825(88)90086-2}

\bibitem[Bends{\o}e and Sigmund(2003)]{Bendsoe2003}
Bends{\o}e MP, Sigmund O (2003) Topology Optimization: Theory, Methods, and Applications. Springer, Berlin Heidelberg. \url{https://doi.org/10.1007/978-3-662-05086-6}

\bibitem[Borrvall and Petersson(2003)]{Borrvall2003}
Borrvall T, Petersson J (2003) Topology optimization of fluids in Stokes flow. Int J Numer Methods Fluids 41:77--107. \url{https://doi.org/10.1002/fld.426}

\bibitem[Choi and Yoon(2024)]{Choi2024}
Choi YH, Yoon GH (2024) A new density filter for pipes for fluid topology optimization. J Fluid Mech 986:A9. \url{https://doi.org/10.1017/jfm.2024.170}

\bibitem[Deng and Chen(2016)]{Deng2016}
Deng J, Chen W (2016) Design for structural flexibility using connected morphable components based topology optimization. Sci China Technol Sci 59:839--851. \url{https://doi.org/10.1007/s11431-016-6027-0}

\bibitem[Dhumal et~al.(2023)]{Dhumal2023}
Dhumal AR, Kulkarni AP, Ambhore NH (2023) A comprehensive review on thermal management of electronic devices. J Eng Appl Sci 70:140. \url{https://doi.org/10.1186/s44147-023-00309-2}

\bibitem[Dilgen et~al.(2018)]{Dilgen2018}
Dilgen CB, Dilgen SB, Fuhrman DR, Sigmund O, Lazarov BS (2018) Topology optimization of turbulent flows. Comput Methods Appl Mech Eng 331:363--393. \url{https://doi.org/10.1016/j.cma.2017.11.029}

\bibitem[Garimella et~al.(2008)]{Garimella2008}
Garimella SV, Fleischer AS, Murthy JY, Keshavarzi A, Prasher R, Patel C, Bhavnani SH, Venkatasubramanian R, Mahajan R, Joshi Y, Sammakia B, Myers BA, Chorosinski L, Baelmans M, Sathyamurthy P, Raad PE (2008) Thermal challenges in next-generation electronic systems. IEEE Trans Components Packag Technol 31(4):801--815. \url{https://doi.org/10.1109/TCAPT.2008.2001197}

\bibitem[Guest et~al.(2004)]{Guest2004}
Guest JK, Pr\'evost JH, Belytschko T (2004) Achieving minimum length scale in topology optimization using nodal design variables and projection functions. Int J Numer Methods Eng 61(2):238--254. \url{https://doi.org/10.1002/nme.1064}

\bibitem[Guo et~al.(2016a)]{Guo2016}
Guo X, Zhang W, Zhang J, Yuan J (2016) Explicit structural topology optimization based on moving morphable components (MMC) with curved skeletons. Comput Methods Appl Mech Eng 310:711--748. \url{https://doi.org/10.1016/j.cma.2016.07.018}

\bibitem[Guo et~al.(2014)]{Guo2014}
Guo X, Zhang W, Zhong W (2014) Doing topology optimization explicitly and geometrically---A new moving morphable components based framework. J Appl Mech 81(8):081009. \url{https://doi.org/10.1115/1.4027609}

\bibitem[Hirotani et~al.(2025)]{Hirotani2025}
Hirotani S, Yaji K, Makihara K, Otsuka K (2025) Data-Driven Real-Time Topology Optimization Using Consistent Rotation-Based Moving Morphable Components. AIAA J 63(10):4491--4497. \url{https://doi.org/10.2514/1.J065458}

\bibitem[Langelaar(2017)]{Langelaar2017}
Langelaar M (2017) An additive manufacturing filter for topology optimization of print-ready designs. Struct Multidisc Optim 55:871--883. \url{https://doi.org/10.1007/s00158-016-1522-2}

\bibitem[Mohammadi and Pironneau(2009)]{Mohammadi2009}
Mohammadi B, Pironneau O (2009) Applied shape optimization for fluids. Oxford University Press, Oxford. \url{https://academic.oup.com/book/1641}

\bibitem[Nguyen et~al.(2010)]{Nguyen2010}
Nguyen TH, Paulino GH, Song J, Le CH (2010) A computational paradigm for multiresolution topology optimization (MTOP). Struct Multidisc Optim 41:525--539. \url{https://doi.org/10.1007/s00158-009-0443-8}

\bibitem[Norato et~al.(2015)]{Norato2015}
Norato JA, Bell BK, Tortorelli DA (2015) A geometry projection method for continuum-based topology optimization with discrete elements. Comput Methods Appl Mech Eng 293:306--327. \url{https://doi.org/10.1016/j.cma.2015.05.005}

\bibitem[Olesen et~al.(2006)]{Olesen2006}
Olesen LH, Okkels F, Bruus H (2006) A high-level programming-language implementation of topology optimization applied to steady-state Navier--Stokes flow. Int J Numer Methods Eng 65:975--1001. \url{https://doi.org/10.1002/nme.1468}

\bibitem[Pejman et~al.(2021)]{Pejman2021}
Pejman R, Sigmund O, Najafi AR (2021) Topology optimization of microvascular composites for active-cooling applications using a geometrical reduced-order model. Struct Multidiscip Optim 64:563--583. \url{https://doi.org/10.1007/s00158-021-02951-x}

\bibitem[Rodrigues et~al.(2002)]{Rodrigues2002}
Rodrigues H, Guedes JM, Bends{\o}e MP (2002) Hierarchical optimization of material and structure. Struct Multidisc Optim 24(1):1--10. \url{https://doi.org/10.1007/s00158-002-0209-z}

\bibitem[Schevenels et~al.(2011)]{Schevenels2011}
Schevenels M, Lazarov BS, Sigmund O (2011) Robust topology optimization accounting for spatially varying manufacturing errors. Comput Methods Appl Mech Eng 200(49--52):3613--3627. \url{https://doi.org/10.1016/j.cma.2011.08.006}

\bibitem[Schury et~al.(2012)]{Schury2012}
Schury F, Stingl M, Wein F (2012) Efficient two-scale optimization of manufacturable graded structures. SIAM J Sci Comput 34(6):B711--B733. \url{https://doi.org/10.1137/110850335}

\bibitem[Sethian and Wiegmann(2000)]{Sethian2000}
Sethian JA, Wiegmann A (2000) Structural boundary design via level set and immersed interface methods. J Comput Phys 163(2):489--528. \url{https://doi.org/10.1006/jcph.2000.6581}

\bibitem[Sigmund(2007)]{Sigmund2007}
Sigmund O (2007) Morphology-based black and white filters for topology optimization. Struct Multidisc Optim 33:401--424. \url{https://doi.org/10.1007/s00158-006-0087-x}

\bibitem[Svanberg(1987)]{Svanberg1987}
Svanberg K (1987) The method of moving asymptotes---a new method for structural optimization. Int J Numer Methods Eng 24(2):359--373. \url{https://doi.org/10.1002/nme.1620240207}

\bibitem[Wang et~al.(2011)]{Wang2011}
Wang F, Lazarov BS, Sigmund O (2011) On projection methods, convergence and robust formulations in topology optimization. Struct Multidisc Optim 43:767--784. \url{https://doi.org/10.1007/s00158-010-0602-y}

\bibitem[Wein et~al.(2020)]{Wein2020}
Wein F, Dunning PD, Norato JA (2020) A review on feature-mapping methods for structural optimization. Struct Multidisc Optim 62(4):1597--1638. \url{https://doi.org/10.1007/s00158-020-02649-6}

\bibitem[Xia and Breitkopf(2014)]{Xia2014}
Xia L, Breitkopf P (2014) Concurrent topology optimization design of material and structure within FE$^{2}$ nonlinear multiscale analysis framework. Comput Methods Appl Mech Eng 278:524--542. \url{https://doi.org/10.1016/j.cma.2014.05.022}

\bibitem[Yan et~al.(2019)]{Yan2019}
Yan S, Wang F, Hong J, Sigmund O (2019) Topology optimization of microchannel heat sinks using a two-layer model. Int J Heat Mass Transfer 143:118462. \url{https://doi.org/10.1016/j.ijheatmasstransfer.2019.118462}

\bibitem[Yoon(2010)]{Yoon2010}
Yoon GH (2010) Topological design of heat dissipating structure with forced convective heat transfer. J Mech Sci Technol 24(6):1225--1233. \url{https://doi.org/10.1007/s12206-010-0328-1}

\bibitem[Yu et~al.(2019)]{Yu2019}
Yu M, Ruan S, Wang X, Li Z, Shen C (2019) Topology optimization of thermal--fluid problem using the MMC-based approach. Struct Multidisc Optim 60:151--165. \url{https://doi.org/10.1007/s00158-019-02206-w}

\bibitem[Zhang et~al.(2016)]{Zhang2016}
Zhang W, Yuan J, Zhang J, Guo X (2016) A new topology optimization approach based on Moving Morphable Components (MMC) and the ersatz material model. Struct Multidisc Optim 53:1243--1260. \url{https://doi.org/10.1007/s00158-015-1372-3}

\end{thebibliography}
\end{document}